\documentclass[11pt]{article}

\usepackage{enumerate}
\usepackage{amsmath,xspace,amssymb,mathrsfs}

\input xy
\xyoption{all}
\xyoption{2cell}
\UseAllTwocells
\CompileMatrices

\title{Homotopy-theoretic aspects of 2-monads}
\author{Stephen Lack%
\thanks{The support of the Australian Research Council and
DETYA is gratefully acknowledged.}
\\School of Computing and Mathematics\\
University of Western Sydney\\
email: {\tt s.lack@uws.edu.au}}
\date{}

\renewcommand{\phi}{\varphi}
\renewcommand{\epsilon}{\varepsilon}

\newcommand{\C}{{\ensuremath{\mathscr C}}\xspace}

\newcommand{\E}{{\ensuremath{\mathscr E}}\xspace}

\newcommand{\G}{{\ensuremath{\mathscr G}}\xspace}

\newcommand{\I}{{\ensuremath{\mathscr I}}\xspace}
\newcommand{\K}{{\ensuremath{\mathscr K}}\xspace}
\newcommand{\LL}{{\ensuremath{\mathscr L}}\xspace}
\newcommand{\M}{{\ensuremath{\mathcal M}}\xspace}

\newcommand{\V}{{\ensuremath{\mathscr V}}\xspace}

\newcommand{\Kf}{\ensuremath{\K_f}\xspace}

\newcommand{\Set}{\textnormal{\bf Set}\xspace}

\newcommand{\SSet}{\textnormal{\bf SSet}\xspace}

\newcommand{\Cat}{\textnormal{\bf Cat}\xspace}

\newcommand{\Catg}{\ensuremath{\Cat_g}\xspace}

\newcommand{\Mndf}{\ensuremath{\textnormal{\bf Mnd}_f(\K)}\xspace}

\newcommand{\Endf}{\ensuremath{\textnormal{\bf End}_f(\K)}\xspace}

\newcommand{\twocat}{\textnormal{\bf 2-Cat}\xspace}
\newcommand{\Twocat}{\textnormal{\bf 2-CAT}\xspace}
\newcommand{\slice}{\Twocat/\K\xspace}
\newcommand{\talg}{{\ensuremath{\textnormal{$T$-Alg}}}\xspace}
\newcommand{\talgs}{{\ensuremath{\textnormal{$T$-Alg}_s}}\xspace}
\newcommand{\pstalg}{{\ensuremath{\textnormal{Ps-$T$-Alg}}}\xspace}
\newcommand{\pstalgs}{{\ensuremath{\textnormal{Ps-$T$-Alg}_s}}\xspace}
\newcommand{\ttalg}{{\ensuremath{\textnormal{$T'$-Alg}}}\xspace}
\newcommand{\ttalgs}{{\ensuremath{\textnormal{$T'$-Alg}_s}}\xspace}
\newcommand{\salg}{{\ensuremath{\textnormal{$S$-Alg}}}\xspace}
\newcommand{\salgs}{{\ensuremath{\textnormal{$S$-Alg}_s}}\xspace}
\newcommand{\malg}{{\ensuremath{\textnormal{$M$-Alg}}}\xspace}
\newcommand{\malgs}{{\ensuremath{\textnormal{$M$-Alg}_s}}\xspace}

\newcommand{\colim}{\textnormal{\sf colim}\xspace}
\newcommand{\Lan}{\textnormal{\sf Lan}\xspace}
\newcommand{\sem}{\textnormal{\sf sem}\xspace}

\newcommand{\id}{\textnormal{\sf id}\xspace}

\newcommand{\Ps}{\textnormal{\sf Ps}\xspace}
\newcommand{\ob}{\textnormal{\sf ob}\xspace}
\newcommand{\iso}{\textnormal{iso}\xspace}

\newcommand{\ct}{\pitchfork}

\newcommand{\bicat}{\textnormal{\bf Bicat}\xspace}

\newcommand{\two}{\ensuremath{{\hbox{\textrm 2}\kern-.25em
        \hbox{\vrule height1.5ex width 0.4pt depth -.2ex}}\kern.2em}\xspace}
\newcommand{\ttwo}{\ensuremath{\two_2}\xspace}

\renewcommand{\>}{\rangle}

\newcommand{\op}{\ensuremath{^{\textnormal{op}}}}

\newcommand{\ot}{\otimes}
\renewcommand{\t}{\times}
\def\endproof{~\hfill$\Box$\vskip 10pt}
\renewcommand{\bar}{\overline}
\newcommand{\Ho}{\textrm{Ho}}

\makeatletter
\newcommand{\pgph}{\paragraph{}}
\renewcommand{\paragraph}{\@startsection
{paragraph}%
{3}%
{0mm}%
{-\baselineskip}%
{-\fontdimen2\font plus -\fontdimen3\font minus -\fontdimen4\font}%
{\normalfont\normalsize\scshape}}%
\makeatother

\numberwithin{equation}{section}
\numberwithin{paragraph}{section}

\newenvironment{theorem}{\paragraph{Theorem:}\em}{\vskip\baselineskip}

\newenvironment{proposition}{\paragraph{Proposition:}\em}{\vskip\baselineskip}
\newenvironment{lemma}{\paragraph{Lemma:}\em}{\vskip\baselineskip}
\newenvironment{remark}{\paragraph{Remark:}}{\vskip\baselineskip}
\newenvironment{example}{\paragraph{Example:}}{\vskip\baselineskip}

\newcommand{\proof}{\noindent{\sc Proof:}\xspace}

\begin{document}

\label{firstpage}
\maketitle

\begin{abstract}
We study 2-monads and their algebras using a \Cat-enriched
version of Quillen model categories, emphasizing the parallels
between the homotopical and 2-categorical points of view. Every
2-category with finite limits and colimits has a canonical model
structure in which the weak equivalences are the equivalences; we
use these to construct more interesting model structures on 2-categories,
including a model structure on the 2-category of algebras for a 2-monad
$T$, and a model structure on a 2-category of 2-monads on a fixed
2-category \K.
\end{abstract}

\section{Introduction}

\pgph
There are obvious connections between 2-category theory and homotopy
theory. It is possible, for instance, to construct a 2-category of
spaces, paths parametrized by intervals of variable length, and
suitably defined equivalence classes of homotopies between them. On
the other hand, in 2-category theory one tends to say that arrows are
isomorphic rather than equal, and that objects are equivalent rather
than isomorphic, typically with some coherence conditions involved,
and this is analogous to working ``up to all higher homotopies''.

In both these cases, the 2-categorical picture is somewhat simpler than
the homotopical one. In the latter case, when one says ``up to all
higher homotopies'', this process does not extend up very far: there
are isomorphisms between arrows, and equations between the isomorphisms,
but that is as far as it goes with 2-categories. For still higher
homotopies, one needs to use $n$-categories for higher $n$, or possibly
$\omega$-categories. This ``degeneracy'' of 2-categories, from the
homotopy point of view, is closely related to the equivalence relations
one must impose on homotopies in order to obtain, as in the previous
paragraph, a 2-category of spaces, paths, and homotopies.

\pgph
Under this analogy, 2-categories correspond to spaces whose homotopy
groups $\pi_n$ are trivial for $n>2$, while mere categories would correspond
to spaces with $\pi_n$ trivial for $n>1$. It is, however, possible to
model all spaces using just categories, via the nerve construction. This
point of view on categories was prominent in
work of Quillen \cite{Quillen-KTheoryI} and Segal \cite{Segal74}, and is the
basis of the theorem of Thomason \cite{Thomason-Cat}
that there is a model structure on the category \Cat
of small categories which is Quillen equivalent to the standard model structure
on the category \SSet of simplicial sets. There is also a corresponding
result for \twocat \cite{Hess-Parent-Tonks-Worytkiewicz}. Under this point
of view, categories are regarded
as being the same if their nerves are homotopy equivalent; this is a
much coarser notion than the usual categorical point of view, adopted
here, that they are the same if they are equivalent as categories.

\pgph
In this paper we pursue the point of view that 2-category theory can
be seen as a slightly degenerate part of homotopy theory, and explore
various consequences, mostly within the area of two-dimensional monad
theory, as developed for example in \cite{BKP}. This point of view builds upon
the earlier papers \cite{qm2cat,qmbicat} which constructed Quillen
model structures on the categories \twocat of 2-categories and 2-functors,
and \bicat of bicategories and strict homomorphisms, investigated their
homotopy-theoretic properties, and related these to the existing theory
of 2-categories.

The difference between this paper and the earlier ones is that
before we looked at a model structure for the category of all
2-categories, whereas here we consider model structures on
particular 2-categories. There is a notion of model structure on an
enriched category, where the base \V for the enrichment itself has a
suitable model structure, and we use this in the case $\V=\Cat$, so
that a \V-category is a 2-category. We therefore speak of a model
\Cat-category, and these are the tools for our analysis of the
homotopy-theoretic aspects of 2-monad theory.

\pgph
A model \Cat-category has three specified classes of morphisms, called
the cofibrations, the weak equivalences, and the fibrations. They satisfy
all the usual properties of model categories, as well as a strengthened
version of the lifting properties, which provides the relationship between
the model structure and the enrichment. We describe the details in
Section~\ref{sect:Cat-model}.

\pgph It turns out that every 2-category with finite limits and
colimits has a canonical model structure in which the weak
equivalences are the equivalences. The details are described in
Section~\ref{sect:trivial}. This can be seen as a
2-categorical analogue of the fact that every category has a model
structure, called the ``trivial'' structure in which the weak
equivalences are the isomorphisms. This trivial structure is almost
never compatible with the enrichment, and so there is little harm in
speaking of the ``trivial model structure on the 2-category \K'' to
mean this one with the equivalences as weak equivalences; for a more
precise statement, see Proposition~\ref{prop:trivial}.  Just as in the
case of ordinary categories, for the trivial model structure on a
2-category, all objects are fibrant and cofibrant. The
factorizations can be constructed in a uniform way using limits and
colimits.

\pgph
The model structures of real interest are not the trivial ones; rather
they can be constructed from the trivial ones via a lifting process. If
\K is a locally presentable 2-category, and $T$ is a 2-monad on \K, there
is a 2-category \talgs of strict $T$-algebras, strict $T$-morphisms, and
$T$-transformations, and a forgetful 2-functor $U_s:\talgs\to\K$ with
a left adjoint $F_s\dashv U_s$. A variety of examples will be discussed
in the following paragraphs; all of these 2-monads and more can be found in the final
section of \cite{BKP}. We use this adjunction to define a model structure
on \talgs, in which a strict $T$-morphism $f$ is a fibration or weak
equivalence in \talgs if and only if the underlying $U_sf$ is one in \K,
where \K is equipped with the trivial model structure. The resulting
model structure on \talgs is not itself trivial; in particular it has
weak equivalences which are not equivalences in \talgs. The details of
the process are developed in Section~\ref{sect:lifted}.

\pgph
An important class of examples is obtained by taking $\K=\Cat$. Then
$T$ describes some sort of algebraic structure borne by categories.
For example this could be monoidal categories, strict monoidal categories,
symmetric monoidal categories, categories with finite limits, categories
with finite products and coproducts satisfying the distributive law,
and so on. In each case, an algebra will consist of a category equipped
with a specific choice of all elements of the structure (for example, a
specific choice of the product $X\times Y$ of two objects, if the structure
involves binary products), and the strict $T$-morphisms will be the
functors which strictly preserve these choices. Such strict morphisms
are of theoretical importance only; usually one would consider the
pseudo $T$-morphisms, which preserve the structure up to suitably coherent
isomorphisms. We write \talg for the 2-category of strict $T$-algebras,
pseudo $T$-morphisms, and $T$-transformations, and usually speak just
of $T$-morphisms, with the ``pseudo'' variety of morphism being the default.
There is a sense, made precise in Theorem~\ref{thm:talg} below, in which
\talg is the ``homotopy 2-category'' of \talgs.

It is familiar in
homotopy theory that up-to-homotopy morphisms from $A$ to $B$ can be
identified, in an up-to-homotopy sense, with ordinary morphisms from a
cofibrant replacement of $A$ to a fibrant replacement of $B$. There is a
corresponding, but rather tighter result here. Every object is already
fibrant, and for any $A$ there is a particular cofibrant replacement $A'$
of $A$ for which the pseudomorphisms from $A$ to $B$ are in bijection
with the strict ones from $A'$ to $B$. In the 2-categorical context the
cofibrant objects are usually called flexible.

\pgph
There are other algebraic structures borne by categories which cannot 
be described in terms of 2-monads on \Cat, but can be described by 2-monads
on the 2-category \Catg of categories, functors, and natural {\em isomorphisms}. 
A typical example is the structure of monoidal closed category. The point
is that the internal hom is covariant in one variable but contravariant in 
the other, and there is no way to describe operations $\C\op\t\C\to\C$ using
2-monads on \Cat. But if we work with \Catg then there is no problem: the 
internal hom is then seen as an operation $\C_\iso\t\C_\iso\to\C_\iso$, where
$\C_\iso$ is the subcategory of \C consisting of the isomorphisms, and
as such the operation is perfectly well-defined. There are subtleties involved, in 
that more work is required to encode the functoriality of the tensor product:
see \cite{BKP}. Similarly such structures as symmetric monoidal closed categories,
compact closed categories, cartesian closed categories, or toposes can be
described by 2-monads on \Catg.

\pgph
Another important class of examples arises on taking $\ob\C$ to be the
objects of a small 2-category \C, and \K to be $[\ob\C,\Cat]$. Then
there is a 2-monad $T$ on \K for which \talgs is the presheaf 2-category
$[\C\op,\Cat]$, consisting of 2-functors, 2-natural transformations,
and modifications. The pseudomorphisms in this case are the pseudonatural
transformations. In Section~\ref{sect:colimits} we study this example,
thinking of presheaves as being the weights for weighted colimits (or
limits). The cofibrant objects are once again important: they are the
weights for flexible colimits. Flexible colimits are those which can be
constructed out of four basic types: coproducts, coinserters, coequifiers,
and splittings of idempotents. The coinserters and coequifiers are
2-dimensional colimits which universally add or make equal 2-cells between
given 1-cells. In the context of \Cat, this corresponds to adding or
making equal morphisms of a category, without changing the objects. One
could use these colimits to force two objects to be isomorphic, but not
to be equal.

It is known that the flexible algebras in \talgs are closed under
flexible colimits; in fact they are precisely the closure under
flexible colimits of the free algebras. In Theorem~\ref{thm:cofibrant}
we give a more general
reason for this first fact: in any model \Cat-category, the cofibrant
objects are always closed under flexible colimits. This is the \Cat-enriched
version of the fact that in any model category the cofibrant objects
are closed under coproducts and retracts.

\pgph
The next example involves a 2-category of 2-monads on a fixed base 2-category
\K, which in this introduction could usefully be taken to be \Cat. The
2-category \Endf of finitary (=filtered-colimit-preserving) endo-2-functors
of \K is locally finitely presentable, and there is a finitary 2-monad $M$
on \Endf whose 2-category \malgs of algebras is the 2-category \Mndf of
finitary 2-monads on \K. In Section~\ref{sect:flexiblemonad} we consider
the trivial model structure on \Endf and the lifted model structure on
\Mndf. One reason for considering \Mndf, is that one can use colimits
in \Mndf to give presentations for monads (exactly as in the unenriched
setting). This depends on the fact that for any object $A$ there is a
2-monad $\<A,A\>$ for which monad morphisms $T\to\<A,A\>$ are in bijection
with $T$-algebra structures on $A$. Thus one can gradually build up
structure on algebras by forming colimits of monads.

The morphisms of \malg are the pseudomorphisms of 2-monads. The main
reason to consider these is to deal with pseudoalgebras. Whereas for
morphisms it is the pseudomorphisms which arise in practice, and the
strict ones are largely just a theoretical construct, it is somewhat
different for algebras. Particular algebraic structures one might
want to consider can most easily be described using strict algebras  ---
for example there is a 2-monad $T$ whose strict algebras are the
not-necessarily-strict monoidal categories --- but some sorts of
formal manipulations one might do fail to preserve strictness, and so
even if one starts with a strict $T$-algebra one might end up with a
non-strict one. This distinction is discussed in
Remark~\ref{rmk:pseudoalgebra}. The connection between pseudomorphisms
of monads and pseudoalgebras, is that to give an object $A$ a
pseudo $T$-algebra structure is equivalent to giving a pseudomorphism
of monads from $T$ to $\<A,A\>$. Using the general technology this
is equivalent to giving a strict map from $T'$ to $\<A,A\>$; that is,
a $T'$-algebra structure on $A$. Furthermore, if the 2-monad $T$ is
flexible, then any pseudo $T$-algebra can be replaced by an {\em isomorphic}
strict one. The fact that flexible colimits of flexible monads are
flexible gives a useful criterion for when a 2-monad given by a
presentation might be flexible.

\pgph The model structure on \Mndf can be used to infer
``semantic'' information about a 2-monad $T$: that is, information
about the 2-category \talg of (strict) $T$-algebras and (pseudo)
$T$-morphisms. The passage from $T$ to \talg is 2-functorial: if
\slice denotes the (enormous!) 2-category of possibly large
2-categories equipped with a 2-functor into \K; the morphisms are
the commutative triangles, then there is a 2-functor
$\sem:\Mndf\op\to\slice$ sending a 2-monad $T$ to \talg, equipped
with the forgetful 2-functor $U:\talg\to\K$; we write
$k^*:\talg\to\salg$ for the 2-functor induced by a morphism of
2-monads $f:S\to T$. Section~\ref{sect:ss} concerns $\sem:\Mndf\op\to\slice$.

The definitions of weak equivalence and fibration for 2-functors coming
from the model structure on \twocat of \cite{qm2cat} make perfectly good
sense for large 2-categories, and it is only size issues which prevent
this from making \slice into a model category, and \sem into a right
Quillen functor. For \sem sends preserves limits, fibrations, and trivial
fibrations (that is, it sends colimits in \Mndf to limits in \slice,
cofibrations in \Mndf to fibrations in \slice, and so on). The extent
to which \sem preserves general weak equivalences is closely related
to the coherence problem for pseudoalgebras (which involves among other
things the replacement of a pseudoalgebra by an equivalent strict one).

The 2-functor $k^*:\talg\to\salg$ restricts to a 2-functor $k^*_s:\talgs\to\salgs$,
which is right adjoint part of a Quillen adjunction, where \talgs and \salgs
are given the lifted model structures. We show in
that $k^*_s$ is a Quillen equivalence if and only if $k^*:\talg\to\salg$
is a biequivalence; that is, if $\sem(k)$ is a weak equivalence in \slice.

\pgph
Instead of monads, another approach to universal algebra
is offered by operads. In \cite{Berger-Moerdijk} operadic analogues
are established to our lifted model structures on \talgs and on \Mndf,
generalizing earlier work by various authors. In one respect, the setting
of \cite{Berger-Moerdijk} is more general than that here, since they
work over an arbitrary monoidal model category \V, whereas here we consider
only the case $\V=\Cat$. This makes for substantial simplifications, due
to the simple nature of the model structure on \Cat. On the other hand,
there are significant simplifications arising from restricting from general
monads to operads, essentially because both in the category of operads, and
in the category of algebras for a given operad, one has much tighter control
over colimits than in the corresponding case for monads. There is also a
more technical difference in that, in contrast to the situation in
\cite{Berger-Moerdijk}, the model structures arising here are not generally
cofibrantly generated, although they are in certain important cases.

In light of this comparison, it is appropriate to give some indication of
the greater generality allowed by monads over operads. Structure described
by operads can only involve operations of the form $A^n\to A$; or, in the
multi-sorted case $A^{n_1}_1\t A^{n_2}_2\t\ldots A^{n_k}_k\to A_m$, where
the superscripts are all natural numbers (corresponding to finite discrete
categories). In the case of monads, one can also use more general limits
such as pullbacks and cotensors. In particular, 2-monads on \Cat allow considering
structures involving maps $A^C\to A$ defined on all diagrams of shape $C$,
for a (not necessarily discrete) category $C$.

\pgph
This paper has had a long gestation period, with the basic results dating back
to 2002. I am grateful to the participants of the seminars at which it was
presented --- the Australian Category Seminar (2002) and the Chicago Category
Seminar (2006) --- for their interest and for various helpful comments. Part
of the writing up was done during a visit to Chicago in May 2006, and I am
very grateful to Peter May, Eugenia Cheng, and the members of the topology/categories
group for their hospitality.

\section{\Cat-model categories}
\label{sect:Cat-model}

\pgph
The category \Cat of small categories and functors has a well-known
``categorical'' or ``folklore'' model structure in which the weak
equivalences are the equivalences of categories, and the fibrations
are the {\em isofibrations}: these are the functors $p:E\to B$ for
which if $e\in E$, and $\beta:b\cong pe$ is an isomorphism in $B$, there
exists an isomorphism $\epsilon:e'\cong e$ in $E$ with $pe'=b$ and
$p\epsilon=\beta$. The model structure is cofibrantly generated,
with generating cofibrations $0\to 1$, $2\to\two$, and $\ttwo\to\two$,
where $2$ is the discrete category with two objects, $\two$ is the
arrow category, and $\ttwo$ is the category with two objects, and a
parallel pair of arrows between them. There is a single generating
trivial cofibration $1\to\I$, where $\I$ is the ``free-living isomorphism''.

\pgph
The cartesian product makes \Cat into a monoidal model category,
in the sense of \cite{Hovey-book}; note that the unit object for the
tensor product is the terminal category $1$, which is cofibrant. We
therefore get, as in \cite{Hovey-book} once again, a notion of model
\Cat-category. Explicitly, a model \Cat-category is a 2-category \K,
with a model structure on the underlying ordinary 2-category $\K_0$
of \K, satisfying the following properties. First of all not just
$\K_0$ but \K must have finite limits and colimits. This reduces to
the further condition that \K have tensors and cotensors by the arrow
category $\two$, which means in turn that for every object $A$ there
are objects $\two\cdot A$ and $\two\ct A$ with natural isomorphisms
$$\K(\two\cdot A,B)\cong\Cat(\two,\K(A,B)\cong\K(A,\two\ct B).$$
As well as this condition on the 2-category, there is also
a compatibility condition on the model structure.
Let $i:A\to B$ be a cofibration and $p:C\to D$
a fibration in \K. Then there is a commutative square
$$\xymatrix{
\K(B,C) \ar[r]^{\K(i,C)} \ar[d]_{\K(B,p)} & \K(A,C) \ar[d]^{\K(A,p)} \\
\K(B,D) \ar[r]_{\K(i,D)} & \K(A,D) }$$
in \Cat, and so an induced functor $[i,p]$ from $\K(B,C)$ to the
pullback. The further property required of a \Cat-model category is that
$[i,p]$ be an isofibration in any case, and moreover an equivalence if either
$i$ or $p$ is trivial.

The fact that $[i,p]$ is surjective on objects if either $i$ or $p$ is
a weak equivalence is just the usual lifting property for the ordinary model
category. We still need (i) that $[i,p]$ is fully faithful if either $i$
or $p$ is a weak equivalence, and (ii) that in any case $[i,p]$ is a
fibration.

Condition (i) says that for any $x,y:B\to C$, there is a bijection between
2-cells $x\to y$ and pairs $\alpha:xi\to yi$ and $\beta:px\to py$ with
$p\alpha=\beta i$. Condition (ii) says that if $z:B\to C$ is given, and
isomorphisms $\alpha:x\cong zi$ and $\beta:y\cong pz$ with $p\alpha=\beta i$,
then there exists a 1-cell $t:B\to E$ and an isomorphism $\sigma:t\cong x$
with $p\sigma=\beta$ and $\sigma i=\alpha$.

The special case $A=0$ of (i) gives the first half of:

\begin{proposition}\label{prop:Cat-model}
If $B$ is cofibrant then $\K(B,-):\K\to\Cat$ preserves fibrations
and trivial fibrations. In particular, if $p:C\to D$ is a trivial
fibration, then composition with $p$ induces an equivalence of
categories $\K(B,C)\simeq\K(B,D)$. Dually, if $E$ is fibrant, then
$\K(-,E):\K\op\to\Cat$ preserves cofibrations and trivial
cofibrations, and if $j:C\to D$ is a trivial cofibration, then
composition with $j$ induces an equivalence $\K(D,E)\simeq\K(C,E)$.
\end{proposition}

\pgph\label{pgph:HoCat}
The homotopy category of \Cat is the category $\Ho\Cat$ of small categories
and isomorphism classes of functors. This category has finite products (computed
as in \Cat), and so we can consider categories enriched in it. Thus the
canonical map $p:\Cat_0\to\Ho\Cat$ preserves finite products, and so every
2-category has an associated $\Ho\Cat$-category, obtained by applying $p$
to each hom-category. (There are corresponding facts with \Cat replaced
by an arbitrary monoidal model category; see \cite{Hovey-book}.)
The homotopy category of a model \Cat-category is canonically a $\Ho\Cat$-category
(once again, see \cite{Hovey-book} for the general situation).
If \K is a model \Cat-category, then the unenriched homotopy
category is the category of objects of \K and isomorphism classes of morphisms.
The enriched homotopy category $\Ho\K$ consists of the objects of \K, and the category
$\Ho(\K(A,B))$ for each $A,B\in\K$. The point is that horizontal composition of
2-cells is now only determined up to isomorphism.

\pgph \label{pgph:Quillen-equivalence}
A right adjoint 2-functor $U:\LL\to\K$ between model \Cat-categories will be called
a right Quillen 2-functor if it sends fibrations to fibrations and
trivial fibrations to trivial fibrations; given that $U$ and not just
its underlying ordinary functor $U_0$ has a left adjoint, this will be
the case if and only if $U_0$ is a right Quillen functor. There is a
derived $\Ho\Cat$-adjunction between the homotopy $\Ho\Cat$-categories,
just as in the usual case. This derived $\Ho\Cat$-adjunction is a
$\Ho\Cat$-equivalence if and only if the unit and counit are invertible,
but this depends only on the underlying ordinary adjunction between
unenriched homotopy categories, so the usual theory of Quillen equivalences
applies.

\section{The trivial \Cat-model structure on a 2-category}
\label{sect:trivial}

\pgph
Let \K be a 2-category with finite limits and colimits. In this section we describe a
\Cat-model structure on \K; we call it the {\em trivial
model structure on the 2-category}. Recall that the trivial model
structure on an ordinary category is obtained by taking the weak
equivalences to be the isomorphisms, and all morphisms to be both
fibrations and cofibrations. This new name can be justified by
Proposition~\ref{prop:trivial} below, which asserts that for a model
\Cat-category \K, if the model structure on the underlying ordinary
category $\K_0$ is trivial then \K is trivial as a model \Cat-category.
(The converse is certainly not true: most trivial model \Cat-categories
will not be trivial at the level of underlying ordinary categories.)

\pgph
The trivial model structure on a 2-category \K can most
concisely be described by saying that a morphism $f:A\to B$
is a weak equivalence or fibration in \K if and only if
the functor $\K(E,f):\K(E,A)\to\K(E,B)$ is one in \Cat, for
every object $E$ of \K. A morphism is a cofibration if and
only if it has the left lifting property with respect to
the trivial fibrations.

Most of this section will be devoted to the
proof of:

\begin{theorem}\label{thm:trivial}
If \K is a 2-category with finite limits and colimits then it becomes
a model \Cat-category with as weak equivalences the (adjoint)
equivalences, and as fibrations the isofibrations. We call this
the {\em trivial model structure}, to distinguish it from any others
which may exist.The factorizations are functorial, and every object
is fibrant and cofibrant.
\end{theorem}

\pgph
First we explicate the definition. A morphism $f:A\to B$ in a 2-category
\K is said to be an equivalence if there exists a morphism $g:B\to A$
with $gf\cong 1_A$ and $fg\cong 1_B$. Since any 2-functor sends equivalences
to equivalences, the equivalences are certainly weak equivalences. Conversely,
if $f:A\to B$ is a weak equivalence, then $\K(B,f):\K(B,A)\to\K(B,B)$ is
an equivalence of categories, so by essential surjectivity there exists a
$g:B\to A$ and $\beta:fg\cong 1_B$. Since $\K(A,f):\K(A,A)\to\K(A,B)$ is also 
an equivalence of categories, and $\K(A,f)gf=fgf\cong f=\K(A,f)1_A$, via
the isomorphism $\beta f$, there is a unique isomorphism $\alpha:gf\cong 1_A$
with $f\alpha=\beta f$. Thus $f$ is an equivalence, and so the weak equivalences
are precisely the equivalences. (We note in passing the well-known fact that
the isomorphisms $gf\cong 1$ and $1\cong fg$ can always be chosen so as to 
satisfy the triangle equations, and so give an adjoint equivalence.)
The weak equivalences are closed under retracts and satisfy the 2-out-of-3 property.

The fibrations are the {\em isofibrations}: these are the
maps $f:A\to B$ such that for any
morphisms $a:X\to A$ and $b:X\to B$, and any invertible 2-cell
$\beta:b\cong fa$, there exists a 1-cell $a':X\to A$ and an invertible
2-cell $\alpha:a'\cong a$ with $fa'=b$ and $f\alpha=\beta$.

It now follows that the trivial fibrations are precisely the retract
equivalences; that is, the morphisms $f:A\to B$ for which there exists
a morphism $g:B\to A$ with $fg=1_A$ and $gf\cong 1_B$. Once again, the
isomorphism can be chosen so as to give an adjoint equivalence.

We define the trivial cofibrations to be the maps with
the left lifting property with respect to the fibrations;
of course these will turn out to be precisely those
cofibrations which are weak equivalences.

\pgph
In the case $\K=\Cat$, this gives the ``categorical'' or ``folklore''
model structure, defined in \cite{Joyal-Tierney-CatE}, for example.
In the case $\K=\Cat^X$ for a set $X$, this gives the pointwise model
structure coming from \Cat.
In the case of $\K=\Cat(\E)$ for a topos \E, this will not in general
be the model structure of \cite{Joyal-Tierney-CatE}, since there the weak
equivalences were the internal functors which are (in the suitably
internal sense) fully faithful and essentially surjective on objects,
and these are more general than the adjoint equivalences unless
the axiom of choice holds in \E. In the case $\K=\Cat(\C)$ for a suitable
finitely complete category \C, it does agree with the model structure
``for the trivial topology'' of \cite{Everaert-Kieboom-VanDerLinden}.

\pgph
The main point of the proof involves a 2-categorical construction
called the pseudolimit of a morphism. If $f:A\to B$ is any 1-cell,
its {\em pseudolimit} is the universal diagram of shape
$$\xymatrix @R1pc {
& A \ar[dd]^{f} \\
L \ar[ur]^{u} \ar[dr]_{v} \rtwocell<\omit>{^\lambda} & {} \\
& B }$$
with $\lambda$ invertible. Thus if $a:X\to A$ and $b:X\to B$ with
$\phi:b\cong fa$, there is a unique 1-cell $c:X\to L$ with $ux=a$,
$vx=b$, and $\lambda x=\phi$. There is also a 2-dimensional aspect
to the universal property, which can most simply be expressed by
saying that if $c,c':X\to L$, then composition with $u$ induces a
bijection between 2-cells $c\to c'$ and 2-cells $uc\to uc'$; in
other words $u$ is representably fully faithful.
Pseudolimits of arrows can be constructed using pullbacks
and cotensors with $\two$.

\begin{proposition}\label{prop:isofibration-pseudolimit}
If an arrow $f:A\to B$ in a 2-category \K admits a pseudolimit
$L$ as above, then $f$ is an isofibration if and only if there
exists a 1-cell $v':L\to A$ and an isomorphism $\lambda':v'\to u$
with $fv'=v$ and $f\lambda'=\lambda$; in other words, if and only
if $\K(L,f)$ is an isofibration in \Cat.
\end{proposition}

\proof The ``only if'' part is immediate; as for the ``if'' part: if
$a:C\to A$, $b:C\to B$, and $\beta:b\cong fa$ are given, let $c:C\to
L$ be the induced map; then $v'c:C\to A$ and $\lambda'c:v'c\cong uc=a$
provide the required lifting.
\endproof

The 1-cells $1:A\to A$ and $f:A\to B$, and the identity 2-cell $f=f$,
induce a unique 1-cell $d:A\to L$ with $ud=1$, $vd=f$, and $\lambda
d=\id_f$. Since $udu=u$, there is a unique invertible 2-cell
$\zeta:du\cong 1$ with $u\zeta$ equal to the identity. It easily
follows that $\zeta d$ is also the identity, so that $d$ and $\zeta$
exhibit $u$ as a retract equivalence.

Furthermore, the universal property of $L$ implies

\begin{proposition}\label{prop:pseudolimit}
The morphism $v:L\to B$ is a fibration.
\end{proposition}

\proof If $c:C\to L$ and $\gamma:e\cong vc$,
then let $(a=uc,b=vc,\beta=\lambda c)$ be the data corresponding
to $c$ via the universal property, so that $\gamma:e\cong
b$. Composing $\beta$ and $\gamma$
gives an isomorphism $e\cong b\cong fa$, which therefore has the
form $\lambda y$ for a unique $y:C\to L$, so that in particular
$vy=e$ and $uy=a$. Now $uy=a=uc$, so there is a unique $\delta:y\cong c$
with $u\delta$ equal to the identity, and now
$$\xymatrix @R1pc {
&& A \ar[dd]^{f} && && A \ar[dd]^{f} && && A \ar[dd]^{f} \\
C \ar[r]^{c} \drrlowertwocell_{e}{^\gamma} & L \ar[ur]^{u} \ar[dr]_{v}
\rtwocell<\omit>{^\lambda} & {} &=&
C \ar[r]^{y} & L \ar[ur]^{u} \ar[dr]_{v} \rtwocell<\omit>{^\lambda} & {} &=&
C \rtwocell^{c}_{y}{^\delta} \drrlowertwocell_{e}{\omit} &
  L \ar[ur]^{u} \ar[dr]_{v}
\rtwocell<\omit>{^\lambda} & {} \\
&& B && && B && && B }$$
and $\lambda c$ is invertible, so that $v\delta=\gamma$, and $\delta$
is the required lifting.
\endproof

\pgph~
\label{para:factorwe}
Observe also that if $f$ is itself an  equivalence, then composing with the
equivalence $u$ gives an equivalence $fu$, so $v$ is an
equivalence since it is isomorphic to $fu$. Thus any morphism $f$ can
be factorized as a weak equivalence $d$ followed by a fibration $v$,
and the fibration will be trivial if (and only if) $f$ was a weak
equivalence.

\pgph
Dually, we can form the pseudocolimit of a morphism $f:A\to B$, involving
$i:A\to C$, $j:B\to C$, and $\lambda:i\cong jf$, as in
$$\xymatrix @R1pc {
A \ar[dd]_{f} \ar[dr]^{i} \\
{}\rtwocell<\omit>{\lambda} & C \\
B \ar[ur]_{j} }$$
and there is an induced
$e:C\to B$ with $ei=f$, $ej=1_B$, $e\lambda=\id$, and an isomorphism
$\epsilon:je\cong 1$ with $\epsilon j$ and $e\epsilon$ both identity 2-cells.
In particular, $e$ is always a trivial fibration.

\begin{lemma}\label{lemma:pseudocolimit}
\begin{enumerate}
\item For any $f$ as above, $i$ is a cofibration and $e$ a
trivial fibration.
\item If $f$ is a weak equivalence, then $i$ is a trivial cofibration.
\end{enumerate}
\end{lemma}

\proof
1. The fact that $e$ is a trivial fibration was observed above. Let's prove
that $i$ is a cofibration. Suppose given then a commutative square
$$\xymatrix{
A \ar[d]_{i} \ar[r]^{u} & E \ar[d]^{p} \\
C \ar[r]_{v} & D }$$
with $p$ a trivial fibration. By the universal property of $C$, to give $v$
is equivalently to give $t:B\to D$ and $\tau:pu\cong tf$.

Since $p$ is a retract equivalence, there exists a 1-cell $s:B\to E$ with
$ps=t$. Then $\tau:pu\cong tf=psf$ has the form $p\sigma$ for a unique
isomorphism $\sigma:u\cong sf$. By the universal property of $C$, there
is a unique $r:C\to E$ with $ri=u$, $rj=s$, and $r\lambda=\sigma$.
If $pr=v$ then $r$ will provide the desired lifting.
But $pri=pu=vi$, $prj=ps=t=vj$, and
$pr\lambda=p\sigma=\tau=v\lambda$, and so $pr=v$ by the universal property
of $C$ once again.

2. We must show that $i$ has the left lifting property with respect
to the fibrations. Suppose given a commutative square
$$\xymatrix{
A \ar[d]_{i} \ar[r]^{u} & E \ar[d]^{p} \\
C \ar[r]_{v} & D }$$
with $p$ a fibration. To give $v$ is equivalently to give
$t:B\to D$ and $\tau:pu\cong tf$. We are assuming that $f$ is
an equivalence, so we can choose $g:B\to A$, $\alpha:1\to gf$,
and $\beta:fg\to 1$ giving an adjoint equivalence. Now there are 1-cells
$ug:B\to E$ and $t:B\to D$, and an isomorphism $t\cong pug$ as in
the left hand side of
$$\xymatrix @R1pc {
&& A \ar[r]^{u} \ar[dd]|{f} & E \ar[dd]^{p} &&
  && A \ar[r]^{u} & E \ar[dd]^{p}\\
& {}\rtwocell<\omit>{\beta} & {}\rtwocell<\omit>{\tau} & {} &=&
  & {}\rtwocell<\omit>{\sigma} & {} \\
B \ar[uurr]^{g} \ar[rr]_{1} && B \ar[r]_{t} & D &&
  B \ar[uurr]^{g} \ar[rr]_{1} \ar@/_1pc/[uurrr]_{s} && B \ar[r]_{t} & D }$$
and so since $p$ is a fibration, there exist an $s:B\to E$ and an
isomorphism $\sigma:ug\cong s$, giving the equality displayed above.
Now by the universal property of $C$, there is a unique $r:C\to E$
with $ri=u$, $rj=s$, and
$$\xymatrix @R1pc {
A \ar[dd]_{f} \ar[dr]^{i} \drruppertwocell^{u}{\omit}  &&&&
A \ar[dd]_{f} \ar[dr]^{1} \drruppertwocell^{u}{\omit} \\
{} \rtwocell<\omit>{\lambda} & C \ar[r]_{r} & E &=&
{} \rtwocell<\omit>{\alpha} & A \ar[r]^{u} & E \\
B \ar[ur]_{j} \urrlowertwocell_{s}{\omit} && &&
B \urrlowertwocell_{s}{\sigma} \ar[ur]_{g} }$$
If $pr=v$ then $r$ will provide the desired lifting.
Now $pri=pu=vi$ and $prj=ps=t=vj$, while
$$\xymatrix @R1pc {
A \ar[dd]_{f} \ar[dr]^{i} && &&
   A \ar[dd]_{f} \ar[dr]^{1} && &&
      A \ar[dd]_{f} \ar[dr]^{1} && && \\
{}\rtwocell<\omit>{\lambda} & C \ar[r]^{r} & E \ar[d]^{p} &=&
   {}\rtwocell<\omit>{\alpha} & A \ar[r]^{u} & E \ar[d]^{p} &=&
      {}\rtwocell<\omit>{\alpha} \drruppertwocell<\omit>{\beta} & A \ar[r]^{u}
      \ar[dr]_{f} \drrtwocell<\omit>{\tau} & E \ar[dr]^{p} \\
B \ar[ur]_{j} && D &&
   B \ar[ur]_{g} \urrlowertwocell_{s}{\sigma} && D &&
      B \ar[rr]_{1}  \ar[ur]_{g} && B \ar[r]_{v} & D }$$
but by one of the triangle equations this last just reduces to
$\tau$; that is to $v\lambda$. Thus $pr\lambda=v\lambda$ and so $pr=v$.
\endproof

\pgph
By the lemma we know that every
morphism can be factorized as a cofibration followed by a trivial fibration,
and that every weak equivalence can be factorized as a trivial cofibration
followed by a trivial fibration. But we saw in \ref{para:factorwe} that every
map can be factorized as a weak equivalence followed by a fibration, thus we
now have both factorization properties. Notice that in order to obtain
the factorization as a trivial cofibration followed by a fibration we
have used both the pseudocolimit and the pseudolimit, whereas for the
other factorization we only needed the pseudocolimit.

\pgph \label{pgph:trivial-cofibration}
We now check that the trivial cofibrations are precisely the maps that
are both weak equivalences and cofibrations.
If $f:A\to B$ is a weak equivalence and a cofibration, then by the lemma
we can factorize it as $f=pi$ with $i$ a trivial cofibration and $p$ a
trivial fibration. The lifting property for cofibrations and trivial fibrations
now makes $f$ a retract of the trivial cofibration $i$, and so $f$ is itself a
trivial cofibration.

If conversely $f$ is a trivial cofibration, then certainly it is a cofibration;
we must show that it is a weak equivalence. To do this, factorize it as a
weak equivalence $d$ followed by a fibration $v$, using the pseudolimit of $f$,
and now the lifting property for trivial cofibrations and fibrations shows
that $f$ is a retract of the weak equivalence $d$, and so a weak equivalence.

This completes the verification of the model category axioms; the
factorizations were explicitly constructed and clearly functorial, and
every object is fibrant and cofibrant. It remains to check that this
gives a model \Cat-category.

\pgph
Let $i:A\to B$ be a cofibration and $p:C\to D$ a fibration in \K.
Let $x,y:B\to C$ be given, with 2-cells $\alpha:xi\to yi$ and
$\beta:px\to py$ satisfying $p\alpha=\beta i$. If $p$ is a trivial
fibration, then in particular it is an equivalence, and so there
is a unique $\gamma:x\to y$ satisfying $p\gamma=\beta$. Furthermore
the 2-cells $\gamma i,\alpha:xi\to yi$ satisfy $p\gamma i=\beta i=p\alpha$,
and so $\gamma i=\alpha$. Similarly if $i$ is an equivalence then there
is a unique $\gamma:x\to y$ satisfying $\gamma i=\alpha$, and it is
also the case that $p\gamma=\beta$. This proves condition (i) for
a model \Cat-category.

As for condition (ii), let $x:A\to C$, $y:B\to D$, $z:B\to C$, $\alpha:x\cong zi$, and
$\beta:y\cong pz$ be given with $p\alpha=\beta i$. We shall verify
the condition first under the assumption that $i$ appears in a
pseudocolimit
$$\xymatrix{
A \ar[rr]^{f} \ar[rd]_{i} & {}\dtwocell<\omit>{\phi} & E \ar[dl]^{j} \\
& B. }$$
The general case will then follow since a general cofibration $i'$ can
be factorized as such an $i$ followed by a trivial fibration, as in
Lemma~\ref{lemma:pseudocolimit}, thus by the lifting property
$i'$ is a retract of $i$, whence the general result.

Suppose then that $i$ does indeed arise in such a pseudocolimit. By
the universal property of $B$, to give $y:B\to D$ is to give
$y_1:A\to D$, $y_2:C\to D$, and an isomorphism $\eta:y_1\cong y_2 f$.
Similarly, to give $z:B\to C$ is to give $z_1:A\to C$, $z_2:E\to C$, and
$\zeta:z_1\cong z_2 f$. To give $\beta:y\cong pz$ is just to give
$\beta_2:y_2\cong pz_2$ (with $\beta_2=\beta j$ and
$\beta i=pz\lambda^{-1}.\beta_2f.y\lambda)$. Thus $\alpha:x\cong zi=z_1$
satisfies $p\alpha=\beta i$ if and only if
$p\alpha=pz\lambda^{-1}.\beta_2f.y\lambda$ or equivalently
$p\zeta.p\alpha=\beta_2f.\eta$.

Since $p$ is a fibration, there exist a 1-cell $t_2:B\to C$ and an
isomorphism $\tau_2:t_2\cong z_2$ as in
$$\xymatrix @R1pc {
& C \ar[dd]^{p} && & C \ar[dd]^{p} \\
&&= \\
B \uurtwocell^{z_2}_{t_2}{^\tau_2} \ar[r]_{y_2} & D  &&
B \ar[r]_{y_2} \uuruppertwocell^{z_2}{^\beta_2} & D }$$
with $pt_2=y_2$ and $p\tau_2=\beta_2$.
The 1-cells $x:A\to C$ and $t_2:B\to C$, and the isomorphism
$$\xymatrix{
x \ar[r]^{\alpha} & z_1 \ar[r]^{\zeta} & z_2 f \ar[r]^{\tau^{-1}_2 f} &
t_2 f }$$
induce a unique $t:B\to E$ by the universal property of the pseudocolimit
$B$. Now $tj=t_2$ and $zj=z_2$, so the isomorphism $\tau_2:t_2\cong z_2$
extends to a unique isomorphism $\tau:t\cong z$ with $\tau j=\tau_2$. 
We shall show that $t$ and $\tau$ have the required properties.

We have $ti=x$ by construction; we show that $pt=y$ using the universal
property of the pseudocolimit $B$. Now $pti=px=yi$ and $ptj=pt_2=y_2=yj$,
while
$$pt\phi=(p\tau^{-1}_2 f).(p\zeta).(p\alpha)=\beta^{-1}_2 f.\beta_2 f.\eta=\eta=y\phi$$
so that $pt=y$ as required.

Thus $p\tau$ and $\beta$ both go from $pt=y$ to $pz$. Since $j$ is
an equivalence, they will be equal if and only if $p\tau j=\beta j$;
but $p\tau j=p\tau_2=\beta_2=\beta j$.

It remains to show that $\tau i=\alpha$. To do so, observe that
$z\phi.\tau i=\tau jf.t\phi$ by the middle-four interchange law,
while $\tau jf.t\phi=\tau_2 f.t\phi=\zeta.\alpha=z\phi.\alpha$,
so that $z\phi.\tau i=z\phi.\alpha$, but $z\phi$ is invertible, so
$\tau i=\alpha$ as desired.

This completes the proof of Theorem~\ref{thm:trivial}.

\begin{proposition}\label{prop:trivial}
Let \K be a model \Cat-category for which the model structure on the
underlying ordinary category is trivial. Then the model \Cat-category
is also trivial.
\end{proposition}

\proof
We shall prove that the only invertible 2-cells are the identities. This
in turn implies that the equivalences are precisely the isomorphisms, and
that all morphisms are isofibrations, and so the proposition will follow.

Let $B$ be an arbitrary object of \K, and let \I be the free-living isomorphism
in \Cat. The cotensor $\I\ct B$ is the ``object of isomorphisms in $B$'', and
is the universal object equipped with morphisms $p,q:\I\ct B\to B$ and an
invertible 2-cell $\theta:p\to q$. There is a unique map $d:B\to\I\ct B$ satisfying
$pd=qd=1$ and $\theta d=\id$. It will suffice to show that $p=q$ 
and $\theta$ is the identity.
Since $q$ is a cofibration,
$\I\ct B$ is fibrant, and $\theta:p\cong 1q$, condition (ii) for a model
\Cat-category implies that there exist a morphism $s:B\to B$ and an isomorphism
$\sigma:s\cong 1$ with $sq=p$ and $\sigma q=\theta$. But
$s=sqd=pd=1$ and $\sigma=\sigma qd=\theta d d=\id$, so $q=p$ and $\theta=\id$
as required.
\endproof


We end with a few degenerate examples.

\begin{example}
If \K is a locally discrete 2-category, meaning that the only 2-cells
are the identities, we may identify it with
its underlying ordinary category. The resulting model structure is
well-known: the weak equivalences
are just the isomorphisms, and all maps are both fibrations and
cofibrations.
\end{example}

\begin{example}
If \K is locally chaotic, so that between any two parallel arrows
$f,g:A\to B$, there is a unique 2-cell $f\to g$ (necessarily invertible),
then once again we can identify \K (in a different way) with its underlying
category. This time a map $f:A\to B$ is a weak equivalence if and only
if there exists an arbitrary map $B\to A$. The trivial fibrations are
precisely the retractions.
\end{example}

\pgph
Finally we observe that even for extremely well-behaved 2-categories, the
trivial model structure need not be cofibrantly generated.
We write $\Cat^{\two}$ for the 2-category of arrows in \Cat: an
object is a functor $a:A\to A'$, a morphism from $a:A\to A'$
to $b:B\to B'$ is a commutative square, involving $f:A\to B$ and
$f':A'\to B'$, and a 2-cell from $(f,f')$ to $(g,g')$ consists
of 2-cells $\alpha:f\to g$ and $\alpha':f'\to g'$ satisfying
$b\alpha=\alpha a'$.

\begin{proposition}
The trivial model structure on $\Cat^{\two}$
is not cofibrantly generated.
\end{proposition}

\proof
Let \K be the category $\Set^{\two}$ of arrows in \Set, seen
as a locally chaotic 2-category. There is a fully faithful 2-functor
$U:\K\to\Cat^{\two}$ sending a function $f:X\to Y$ to the corresponding
functor between the chaotic categories on $X$ and $Y$. This 2-functor
has a left adjoint, which sends a functor to the corresponding function
between object-sets. The trivial model structure on $\Set^{\two}$ can
be obtained as the lifting of the trivial model structure on $\Cat^{\two}$.
Now if $\Cat^{\two}$ were cofibrantly generated, then so would be
$\Set^{\two}$, but it is not. For if $(m_i:A_i\to A_i+B_i)_{i\in I}$
were a small family of generating cofibrations (here $A_i$ and $B_i$
denote objects of $\Set^{\two}$; that is, functions) then so would be
$(0\to B_i)_{i\in I}$. Now every object is cofibrant, so would
be a retract of a coproduct of $B_i$'s. By the exactness properties
of $\Set$, it would then follow that every object was a coproduct
of retracts of $B_i$'s. Now the closure under retracts of the
$B_i$'s is still small, so there would be a small full subcategory
\G of $\Set^{\two}$ such that every object was a coproduct of
objects in \G. But now consider the objects of the form $X\to 1$.
These constitute a large family, and none of them can be decomposed
non-trivially as a coproduct. Thus there cannot be a small family
of generating cofibrations.
\endproof

This should not perhaps be too surprising. We have defined the weak
equivalences and fibrations by lifting through the representables
$\K(C,-)$ for arbitrary $C$. Since this is generally a large set of
objects it is not surprising that it would not lead to a cofibrantly
generated structure. In some cases however a small set of objects
will suffice, and then the structure will be cofibrantly generated.
In particular, if $\K=\Cat$, then it suffices to use just the single
representable $\Cat(1,-)$ (which is the identity 2-functor on \Cat).

\section{The lifted model structure on a 2-category of algebras}
\label{sect:lifted}

\pgph
Suppose now that \K is a locally presentable 2-category,
endowed with the trivial model structure as in the previous section.
Suppose that $T$ is a 2-monad on \K with rank (preserves $\alpha$-filtered
colimits for some $\alpha$), and that $\talgs$ is
the 2-category of strict $T$-algebras, strict morphisms,
and $T$-transformations. Then the forgetful 2-functor
$U_s:\talgs\to\K$ has a left 2-adjoint $F_s$, with unit $n:1\to U_sF_s$
and counit $e:F_sU_s\to 1$. We shall use this adjunction to construct
a ``lifted'' model structure on \talgs. A morphism $f$ in \talgs is
defined to be a fibration or weak equivalence if and only if $U_s f$
is one in \K, while a morphism is a cofibration if and only if it has the left
lifting property with respect to the trivial fibrations (the maps
which are both fibrations and weak equivalences), and a trivial
cofibration if and only if it has the left lifting property with
respect to the fibrations. The fact that this is a model \Cat-structure
will follow immediately from the fact that it is a model structure,
thanks to the 2-dimensional aspect of the 2-adjunction.

There exist many theorems about lifting model structures, but they
generally depend upon the lifted model structure being cofibrantly
generated, which we are not assuming here.

\pgph \label{pgph:we-fibration}
Given a strict morphism $f:A\to B$ we can form in \talgs the pseudolimit
$L$ of $f$, with projections $u:L\to A$ and
$v:L\to B$ and isomorphism $\lambda:v\cong fu$; since
$U_s:\talgs\to\K$ preserves limits, $U_sL$ will also be the limit in \K,
and so $U_s v:U_sL\to U_sB$ will be a fibration in \K by
Proposition~\ref{prop:pseudolimit}, and so in turn $v$ will
be a fibration in \talgs. Furthermore, we have the unique induced
$d:A\to L$ with $ud=1$ and $\phi d$ the identity, and just as
before this $d$ is an equivalence in \talgs, and so in particular a
weak equivalence. This proves that every map can be factorized as
a (weak) equivalence followed by fibration. This already implies that
every trivial cofibration is a weak equivalence, by the same argument
used in Section~\ref{pgph:trivial-cofibration}.

\pgph
Let $f:A\to B$ be an arbitrary strict morphism. Factorize $U_sf:U_sA\to U_sB$
in \K as a cofibration $i_1:U_sA\to X_1$ followed by a trivial fibration
$p_1:X_1\to U_sB$. Pushout $F_si_1$ along the counit $eA:F_sU_sA\to A$, and
form the induced map $f_1$ as in
$$\xymatrix{
F_sU_sA \ar[r]^{F_si_1} \ar[d]_{eA} & F_sX_1 \ar[d]^{c_1} \ar[r]^{F_sp_1} &
F_sU_sB \ar[d]^{eB} \\
A \ar[r]_{j_1} & C_1 \ar[r]_{f_1} & B}$$
where the left square is the pushout, and $f_1j_1=f$. Now $i_1$ is
a trivial cofibration in \K, so $F_si_1$ is a trivial cofibration in
\talgs, and so in turn is its pushout $j_1$. There is not so much
we can say about $f_1$ at this stage, but we do know that $U_sf_1$ has
a section, for if $s_1$ is a section of the trivial fibration $p_1$,
then
$U_sf_1.U_sc_1.U_sF_ss_1.nU_sB=U_seB.U_sF_sp_1.U_sF_ss_1.nU_sB=U_seB.nU_sB=1$.

If $f$ is in fact a weak equivalence, then since $f=f_1 j_1$ and $j_1$
is a weak equivalence, $f_1$ will be one too. But now $U_sf_1$ is
a weak equivalence with a section, hence a trivial fibration, and
so finally $f_1$ is a trivial fibration. Thus every weak equivalence
factorizes as trivial cofibration followed by a trivial fibration.
Combined with the factorization, given in Section~\ref{pgph:we-fibration},
of any map into a weak equivalence followed by a fibration, this now
proves that any map can be factorized as a trivial cofibration followed
by a fibration. It now follows that the trivial cofibrations are
precisely the weak equivalences which are cofibrations, by the argument
used in Section~\ref{pgph:trivial-cofibration}.

\pgph
So far things have gone essentially as usual. It remains to show
the existence of the other factorization: cofibration followed by
trivial fibration. This is the most technical part of the proof.
Suppose again then that $f$ is arbitrary, and
factorize it as $f=f_1j_1$ as above. We know that $U_sf_1$ has a
section. If we could show that for any two maps $x,y:U_sB\to C_1$
with $U_sf_1.x\cong U_sf_1.y$ we
have $x\cong y$, then $U_sf_1$ would be a trivial fibration, and we
would be done; but in general there is no reason why this should
be true, and there is more work to be done. Factorize $f_1$ as $f_2j_2$
via the same process, and now iterate to obtain a trivial cofibration
$j_{n+1}:C_n\to C_{n+1}$ and map $f_{n+1}:C_{n+1}\to B$ for any $n$. Continue
transfinitely, setting $C_m=\colim_{n<m}C_n$ for any limit ordinal $m$.
Any transfinite composite of the $j$'s will be a trivial cofibration,
and each $f_n$ will have the property that $U_sf_n$ has a section.
So if we can find an $n$ such that for any $x,y:U_sB\to C_n$ with
$U_sf_n.x\cong U_sf_n.y$ we have $x\cong y$, then $f_n$ will be a trivial
fibration, and we will be done.

Let $\alpha$ be a regular cardinal with the property that $T$
preserves $\alpha$-filtered colimits and that $U_sB$ is
$\alpha$-presentable in \K. Since $T$ preserves $\alpha$-filtered
colimits, so does $U_s$. Let $x,y:U_sB\to C_{\alpha}$ be given with
$U_sf_{\alpha}.x\cong U_sf_{\alpha}.y$. Now
$C_{\alpha}=\colim_{n<\alpha} C_i$, so there exists an $n<\alpha$
such that $x$ and $y$ land in $C_n$, say via $x',y':U_sB\to U_sC_n$.
Now $U_sf_n.x'=U_sf_{\alpha}.x\cong U_sf_{\alpha}.y=U_sf_n.y'$, and
$U_sf_n=p_{n+1}.i_{n+1}$, with $p_{n+1}$ a trivial fibration, so the
isomorphism lifts through $p_{n+1}$ to give $i_{n+1}.x'\cong
i_{n+1}.y'$. Now
$$j_{n+1}.eC_n.F_sx'=c_{n+1}.F_si_{n+1}.F_sx'\cong c_{n+1}.F_si_{n+1}.F_sy'=
j_{n+1}.eC_n.F_sy'$$
and $j_{n+1}$ is a trivial cofibration, so has a retraction, and
so $eC_n.F_sx'\cong eC_n.F_sy'$, but by adjointness this is just
$x'\cong y'$, which finally gives $x\cong y$ as required.

This proves the existence of the model structure; it is automatically
a model \Cat-structure, via the 2-dimensional aspect of the adjunction.

\begin{theorem}
For a 2-monad $T$ with rank, on a locally finitely presentable 2-category
\K, the category \talgs of strict $T$-algebras and strict $T$-morphisms
has a cofibrantly generated \Cat-model structure for which the weak
equivalences are the maps which are equivalences in \K, and the fibrations
are the maps which are isofibrations in \K.
\end{theorem}

\begin{remark}
Observe that the class of all strict $T$-morphisms of the form
$Ti:TX\to TZ$ with $i:X\to Z$ a cofibration in \K, while not small,
does nonetheless generate the cofibrations of \talgs. Furthermore,
we can even restrict to those $Ti$ for which there exist $f:X\to Y$,
$j:Y\to Z$ , and $\lambda:i\cong fj$, such that $i$, $j$, and $\lambda$
exhibit $Z$ is as the pseudocolimit in \K of $f$, for every cofibration
in \K is a retract of one of these.
\end{remark}

Once again every object is fibrant, but it is no longer the case that
every object is cofibrant; we shall see below that the cofibrant objects
are precisely the {\em flexible} ones, in the sense of \cite{BKP}.

\pgph
The strict morphisms, as in \talgs, are very useful for theoretical
reasons, but in practice they are rare. More common are the {\em pseudo}
$T$-morphisms, which preserve the structure only up to coherent isomorphism.
Since we are treating this as the basic notion of morphism, we call them
simply $T$-morphisms. They are the morphisms of a 2-category \talg of
(still strict) $T$-algebras, $T$-morphisms, and $T$-transformations. The
inclusion 2-functor $\talgs\to\talg$ is the identity on objects, and
fully faithful on the hom-categories.

\begin{remark}\label{rmk:pseudoalgebra}
There is also a notion of {\em pseudo $T$-algebra} for a 2-monad $T$,
in which the usual laws for $T$-algebras are replaced by coherent
isomorphisms; these are considered in Section~\ref{sect:flexiblemonad}
below. The pseudo algebras are less important than the strict ones
for two reasons. First there is a ``theoretical'' reason, discussed
in Section~\ref{sect:flexiblemonad}, that for a 2-monad $T$ with rank
on a locally presentable 2-category \K, the pseudo $T$-algebras are
just the strict algebras for a different 2-monad $T'$. There is also
a more practical reason, which we illustrate with the example of monoidal
categories. There is a 2-monad $T$ on \Cat whose strict algebras are
the strict monoidal categories. It is true that ``up to equivalence''
the pseudo $T$-algebras are the same as the (not necessarily strict)
monoidal categories, but this is a relatively hard fact. Much easier
is the fact that there is a 2-monad $S$ whose strict algebras are
precisely the monoidal categories (we sketch below the slightly simpler
case of ``semigroupoidal categories''). The situation is similar but
more pronounced with more complicated structures than monoidal categories.
The reason for considering pseudoalgebras at all is that some natural constructions
on algebras only produce pseudoalgebras, even if one starts with a strict
one.
\end{remark}

\pgph
We write $U:\talg\to\K$ for the forgetful map, to
distinguish it from $U_s:\talgs\to\K$. The evident inclusion $J:\talgs\to\talg$
clearly satisfies $UJ=U_s$. It was proved in \cite{BKP} that
$J:\talgs\to\talg$ has a left 2-adjoint, sending an algebra $A$ to an
algebra $A'$, with counit a strict map $q:A'\to A$ and with unit a
pseudomap $p:A\to A'$ satisfying $qp=1$ and $pq\cong 1$. Thus $q$ is
a trivial fibration. The universal
property of $A'$ asserts among other things that for any algebra $B$,
composition with $p$ induces a bijection between strict maps $A'\to B$
and pseudo maps $A\to B$.

\begin{proposition}\label{prop:equivalence}
For a strict morphism $f:A\to B$, the following are equivalent:
\begin{enumerate}[(i)]
\item $f$ is a weak equivalence;
\item $U_s f$ is an equivalence in \K;
\item $J f$ is an equivalence in \talg.
\end{enumerate}
\end{proposition}

\proof The equivalence of (i) and (ii) holds by definition of
weak equivalences in \talgs. The equivalence of (ii) and (iii)
is a routine (but important) exercise.
\endproof

\begin{proposition}
For a strict morphism $f:A\to B$, the following are equivalent:
\begin{enumerate}[(i)]
\item $f$ is a fibration;
\item $U_s f$ is an isofibration in \K;
\item $J f$ is an isofibration in \talg.
\end{enumerate}
\end{proposition}

\proof The equivalence of (i) and (ii) holds by definition of
weak equivalences in \talgs. The equivalence of (ii) and (iii)
is a straightforward consequence of
Proposition~\ref{prop:isofibration-pseudolimit}.
\endproof

Notice in particular that $q:A'\to A$ is a trivial fibration for
any algebra $A$, since we have $qp=1$ and $pq\cong 1$ in \talg (and
in \K). Thus for any algebra $A$, the map $q:A'\to A$ has a section in
\talg; if it has a section in \talgs --- that is, a {\em strict} map
$r:A\to A'$ with $qr=1$ --- then $A$ is said to be {\em flexible}
\cite{BKP}.

\begin{theorem}
The cofibrant objects of \talgs are precisely the flexible algebras;
in particular, any algebra of the form $A'$ is cofibrant, and so a
cofibrant replacement for $A$. Every free algebra is flexible.
\end{theorem}

\proof Since $q:A'\to A$ is a trivial fibration, then certainly it
will have a section if $A$ is cofibrant. Thus cofibrant objects are
flexible. For the converse, it will suffice to show that each $A'$
is cofibrant, for any retract of a cofibrant object is cofibrant.

Suppose then that $t:E\to B$ is a trivial fibration in \talgs, and
$v:A'\to B$ an arbitrary strict map. We must show that it lifts through
$t$. There is a pseudomorphism $s:B\to E$ with $ts=1$, and so $tsv=v$.
But the pseudomorphism $svp:A\to E$ has the form $up$ for a unique
strict map $u:A'\to E$. Now the strict maps $tu$ and $v$ from $A'$ to $B$
satisfy $tup=tsvp=vp$, and so $tu=v$ by the universal property of $A'$,
which gives the required lifting.

To see that free algebras are flexible, observe that any object
$X\in\K$ is cofibrant, but for the lifted model structure the left
adjoint preserves cofibrant objects, so the free algebra $TX$ on $X$
is cofibrant, and so flexible.
\endproof

The same relationship between flexibility and cofibrancy was observed
in \cite{qm2cat}.

\begin{proposition}\label{prop:flexibledomain}
Any pseudomorphism with flexible domain is isomorphic to a strict
morphism.
\end{proposition}

\proof
Let $r:A\to A'$ be a strict morphism which is a section of $q:A'\to A$.
Now $qrq=q=q1$, and $q$ is an equivalence in \talg, so $rq\cong 1$ in
\talg; but $rq$ and $1$ are in \talgs, and the inclusion $J:\talgs\to\talg$
is locally fully faithful (fully faithful on 2-cells) and so $rq\cong 1$
in \talgs. This also implies that $r=rqp\cong p$ in \talg.

Now suppose that $f:A\to B$ is a pseudomorphism. It can be written
as $f=gp$ for a unique strict morphism $g:A'\to B$, and now
$gr\cong gp=f$, so that the pseudomorphism $f$ is isomorphic to the strict
morphism $gr$.
\endproof

\pgph
The homotopy category of \talgs is easily described: it is
the category of strict $T$-algebras, and isomorphism classes of pseudo
$T$-morphisms. As explained in Section~\ref{pgph:HoCat}, it has a canonical
enrichment to a $\Ho\Cat$-category. But the resulting $\Ho\Cat$-category
can also be seen as the $\Ho\Cat$-category underlying the 2-category \talg (again
in the sense of Section~\ref{pgph:HoCat}). Thus \talg is a kind of
``homotopy 2-category'' of \talgs. This point of view is reinforced
by the following proposition, which describes a universal property of \talg.

For 2-categories \M and \LL we write $\Ps(\M,\LL)$ for the 2-category of 2-functors,
pseudonatural transformations, and modifications, from \M to \LL.

\begin{theorem}\label{thm:talg}
Let \LL be any 2-category. Composition with $J:\talgs\to\talg$ induces
a biequivalence of 2-categories between $\Ps(\talg,\LL)$ and the
full sub-2-category $\Ps_w(\talgs,\LL)$ consisting of those
2-functors $\talgs\to\LL$ sending weak equivalences to equivalences.
\end{theorem}

\proof
First observe $J:\talgs\to\talg$ sends weak equivalences to equivalences,
by Proposition~\ref{prop:equivalence}, and pseudofunctors preserve
equivalences, so composition with $J$
does indeed induce a 2-functor $R:\Ps(\talg,\LL)\to\Ps_w(\talgs,\LL)$.
If $F:\talgs\to\LL$ sends weak equivalences to equivalences, then in
particular it sends each $q_A:A'\to A$ to a weak equivalence. Now
$q$ is the counit of the adjunction $L\dashv J$, so $Fq$
is an equivalence in $\Ps(\talgs,\LL)$, and so $F\simeq FLJ=R(FL)$,
and $R$ is biessentially surjective on objects. To see that it is
an equivalence on hom-categories, and so a biequivalence, observe
that for  $M,N:\talg\to\LL$ we have
\begin{align*}
\Ps_w(\talgs,\LL)(MJ,NJ) &= \Ps(\talgs,\LL)(MJ,NJ) \\
                           &\simeq \Ps(\talg,\LL)(M,NJL) \\
                           &\simeq \Ps(\talg,\LL)(M,N)
\end{align*}
using adjointness and the fact that the unit $1\to JL$ is an equivalence.
\endproof

We cannot expect this to work using 2-natural transformations. It is clear
from the proof that rather than sending all weak equivalences to equivalences,
we could ask only that trivial fibrations be sent to equivalences; but by
Ken Brown's lemma \cite[1.1.2]{Hovey-book} and the fact that all objects
of \talgs are fibrant, any 2-functor $\talgs\to\LL$ sending all trivial fibrations
to equivalences must in fact send all weak equivalences to equivalences.

\section{Flexible colimits}
\label{sect:colimits}

\pgph
In this section we consider \talgs and its model structure for a particular
case of $T$, relevant to (weighted) colimits in 2-categories. Recall that if
$S:\C\to\K$ and $J:\C\op\to\Cat$ are 2-functors, with \C small, we
write $J*S$ for the $J$-weighted colimit of $S$, defined by an isomorphism
$$\K(J*S,A)\cong[\C\op,\Cat](J,\K(S,A))$$
natural in $A$, where $[\C\op,\Cat]$ is the 2-category of 2-functors,
2-natural transformations, and modifications. The presheaf $J$ is
called the weight.
We shall describe a 2-category \K and a 2-monad $T$ on \K for which
\talgs is precisely this 2-category $[\C\op,\Cat]$.

\pgph
Let \C be a small 2-category, and write $\ob\C$ for its set
of objects. Our base 2-category \K will be $[\ob\C,\Cat]$; this is just
the set of $\ob\C$-indexed families of categories. The 2-category
$[\C\op,\Cat]$ has an evident forgetful 2-functor
$U_s:[\C\op,\Cat]\to[\ob\C,\Cat]$ given by restriction along the inclusion
$\ob\C\to\C\op$, and $U_s$ has left and right adjoints given by left and right
Kan extension along the inclusion. It is now straightforward to verify using
(an enriched variant of) Beck's theorem that $U_s$ is monadic, so that
$[\C\op,\Cat]$ has the form \talgs for a 2-monad $T$ on $[\ob\C,\Cat]$.
Since $U_s$ has a right as well as a left adjoint, it follows that $T$
too has a right adjoint, so preserves all colimits, thus in particular
is finitary. The
corresponding 2-category \talg is $\Ps(\C\op,\Cat)$, consisting of the
2-functors, pseudonatural transformations, and modifications. A flexible
algebra for this 2-monad is called a {\em flexible weight} \cite{BKPS},
and colimits weighted by flexible weights are called {\em flexible colimits}.

Just as all ordinary colimits can be computed using coproducts and
coequalizers, all flexible colimits can be computed using four basic
types of flexible colimit: coproducts,
coinserters, coequifiers, and splittings of idempotents \cite{BKPS}.
For a good introduction to these various limit notions, see
\cite{Kelly-limits}. Here we simply recall that all these notions are
defined by universal properties involving 2-natural isomorphisms, and
that the coinserter of a pair $f,g:A\to B$ is the universal map $p:B\to C$
equipped with a 2-cell $pf\to pg$; while the coequifier of a parallel
pair $\alpha$ and $\beta$ of 2-cells between 1-cells $f,g:A\to B$ is the
universal $p:B\to C$ for which $p\alpha=p\beta$.

\pgph
A fundamental result is that the flexible algebras are closed under
flexible colimits (in \talgs). This is equivalent to being closed
under these four types of colimit. Here we offer an alternative
viewpoint on this fundamental result, based on the fact that the
flexible algebras are precisely the cofibrant objects in a model
\Cat-category. In any model category the cofibrant objects are
closed under coproducts and splittings of idempotents; but in a
model \Cat-category we have:

\begin{theorem}\label{thm:cofibrant}
In a model \Cat-category, the cofibrant objects are closed under
flexible colimits.
\end{theorem}

\proof
Since cofibrant objects are always closed under
coproducts and splittings of idempotents (retracts), it remains to
show that they are closed under coinserters and coequifiers.

Let $f,g:F\to A$ be morphisms between cofibrant objects, and let $i:A\to B$
and $\alpha:if\to ig$ exhibit $B$ as the coinserter of $f$ and $g$.
We shall show that $i$ is a cofibration, and so that $B$ is cofibrant.
Suppose then that
$$\xymatrix{
A \ar[r]^{x} \ar[d]_{i} & C \ar[d]^{p} \\
B \ar[r]_{y} & D }$$
is a commutative square with $p$ a trivial fibration. Regarding $p$ and
$x$ as fixed, to give $y$ is just to give $\phi:pxf\to pxg$.
Since $A$ is cofibrant there is  by Proposition~\ref{prop:Cat-model}
a unique 2-cell $\psi:xf\to xg$ with $p\psi=\phi$.
By the universal property of the coinserter there is now a unique
$z:B\to C$ with $zi=x$ and $z\alpha=\psi$. On the other hand
$pzi=px=yi$ and $pz\alpha=p\psi=\phi=y\alpha$ so $pz=y$ by the
uniqueness part of the universal property. Thus $z$ is the required
fill-in and so $i$ is a cofibration.

Now we turn to coequifiers. Let $f,g:F\to A$ be morphisms between
cofibrant objects, and let $i:A\to B$ be the coequifier of 2-cells
$\alpha,\beta:f\to g$. We shall show that $i$ is a cofibration and
so that $B$ is cofibrant. Consider a square as above; this time $y$
is uniquely determined by $px$, and its existence is equivalent to
the equation $px\alpha=px\beta$. Since $F$ is cofibrant and $p$ is
a trivial fibration, we have $x\alpha=x\beta$ by
Proposition~\ref{prop:Cat-model} once again, and so a unique
$z:B\to C$ with $zi=x$. The equation $pz=y$ is immediate consequence
using the universal property of $B$ once again, so $z$ provides the
fill-in for the square, and $i$ is once again a cofibration.
\endproof

\begin{remark}
In the case of the lifted model structure on \talgs, the flexible
(=cofibrant) algebras are precisely the closure under flexible colimits
of the free algebras. On the one hand, all objects of \K are cofibrant,
so all free algebras are cofibrant, and we have seen that cofibrant
objects are closed under flexible colimits. For the converse, it
was shown in \cite{codescent} that for any algebra $A$, the algebra
$A'$ can be constructed from free algebras using coinserters and
coequifiers (in \talgs); since the flexible algebras are the retracts
of the $A'$, it follows that they are flexible colimits of free algebras.
\end{remark}

\section{Flexible monads}
\label{sect:flexiblemonad}

\pgph
In this section we study a certain 2-category of 2-monads, and a lifted
model structure coming from the underlying 2-category of endo-2-functors.
We continue to consider a fixed locally finitely presentable
2-category \K; for this section the most important case is
$\K=\Cat$. A 2-functor $T:\K\to\K$ is said to be {\em finitary} if
it preserves filtered colimits; or, equivalently, if it is the left
Kan extension of its restriction to the full sub-2-category \Kf of
\K consisting of the finitely presentable objects. Since the
composite of finitary 2-functors is clearly finitary, and identity
2-functors are so, one obtains a strict monoidal category \Endf of
finitary endo-2-functors of \K. As a category, it is equivalent to
the 2-functor category $[\Kf,\K]$: the equivalence sends a finitary
2-functor $T$ to its composite $TJ$ with the inclusion $J:\Kf\to\K$,
and sends $S:\Kf\to\K$ to the left Kan extension $\Lan_J(S)$. The
strict monoidal structure on \Endf transports across the equivalence
to give a (no longer strict) monoidal structure on $[\Kf,\K]$: the
tensor product $S\circ R$ is given by $\Lan_J(S)R$, and the unit is
$J$. We sometimes identify 2-functors $\Kf\to\K$ with the
corresponding finitary endo-2-functors of \K.

\pgph
A 2-monad on \K consists of a 2-functor $T:\K\to\K$ equipped
with 2-natural transformations $m:T^2\to T$ and $i:1\to T$
satisfying the usual monad equations. It is said to be finitary
if the endo-2-functor $T$ is so. (We often allow ourselves to
speak of ``the 2-monad $T$'', leaving the multiplication $m$ and
unit $i$ understood.)

There is now a category \Mndf of finitary 2-monads on \K and
{\em strict} morphisms; it is
the category of monoids in \Endf or equivalently in $[\Kf,\K]$).
The forgetful functor $W:\Mndf\to[\Kf,\K]$ has a left adjoint
$H$ and is monadic.

\pgph \label{pgph:<A,A>}
Monad morphisms are useful for describing
algebras for monads. Recall from \cite{Kelly-LaxAlg} or
\cite[Section~2]{property} that if $A$ and $B$ are objects of \K
then there is a 2-functor $\<A,B\>:\Kf\to\K$ which sends a finitely
presentable object $C$ to the cotensor $\K(C,A)\ct B$, and now to
give a 2-natural transformation $T\to\<A,B\>$ is equivalently to
give a map $TA\to B$ in \K; more precisely, we have an isomorphism
of categories $[\Kf,\K](T,\<A,B\>)\cong\K(TA,B)$. Furthermore, if
$A=B$, then there is a 2-monad structure on $\<A,A\>$, such that for
a 2-monad $T$, a 2-natural transformation $T\to\<A,A\>$ is a monad
map if and only if the corresponding $TA\to A$ makes $A$ into a
$T$-algebra. This observation illustrates the importance of colimits
in \Mndf: it shows, for example, that an algebra for the coproduct
$S+T$ of monads $S$ and $T$ is just an object equipped with an
$S$-algebra structure and a $T$-algebra structure. We shall see
further examples below.

First observe, following \cite{Kelly-LaxAlg,property} once again,
that if $f,g:A\to B$ we may form the comma-object
$$\xymatrix{
\{f,g\} \ar[r]^{c} \ar[d]_{d} \drtwocell<\omit>{\lambda} &
\<A,A\> \ar[d]^{\<A,g\>} \\
\<B,B\> \ar[r]_{\<f,B\>} & \<A,B\> }$$ in \Endf, and now to give a
2-natural transformation $\gamma:T\to\{f,g\}$ is equivalently to
give morphisms $a:TA\to A$ (corresponding to $c\gamma$) and $b:TB\to
B$ (corresponding to $d\gamma$), and an invertible 2-cell $b.Tf\to
ga$. Once again, if $f=g$, then there is a trivial 2-monad structure
on $\{f,f\}$ such that if $T$ is a 2-monad, then $\gamma$ is a monad
map if and only if $(A,a)$ and $(B,b)$ are $T$-algebras and the
induced 2-cell $\bar{f}:b.Tf\to fa$ makes $(f,\bar{f})$ into a
$T$-morphism. Thus once again we can analyze the (pseudo)morphisms
of algebras for colimits of 2-monads. Finally, if $\rho:f\to g$ is a
2-cell in \K, then we may form the pullback
$$\xymatrix{
[\rho,\rho] \ar[r] \ar[d] & \{f,f\} \ar[d]^{\{f,\rho\}} \\
\{g,g\} \ar[r]_{\{\rho,g\}} & \{f,g\} }$$ and this has a canonical
monad structure for which monad maps $T\to[\rho,\rho]$ correspond to
2-cells in \talg.

\pgph
This allows us to give {\em presentations} for 2-monads, as in \cite{KP}.
For example, take $\K=\Cat$, and let $E:\Cat\to\Cat$ be the finitary
2-functor sending a category $\C$ to $\C\t\C$. An algebra for the
free monad $HE$ on $E$ is just a category \C with a functor $\ot:\C\t\C\to\C$.
A (pseudo)morphism is a functor between such categories which
preserves the ``tensor product'' up to an arbitrary natural
isomorphism: there are no coherence conditions at this stage.

Now let $D:\Cat\to\Cat$ be the finitary 2-functor sending $\C$ to
$\C\t\C\t\C$, and $HD$ the free monad on $D$. Since for any
$HE$-algebra \C there are two trivial $HD$-algebra structures
(corresponding to the two bracketings of a triple product), there
are two induced monad maps $f,g:HD\to HE$. We can form the {\em
co-isoinserter} of these maps, which is the universal monad map
$r:HE\to S$ equipped with a monad isomorphism $\rho:rf\cong rg$. An
$S$-algebra is now a category \C, with a functor $\ot:\C\t\C\to\C$,
and a natural isomorphism $\alpha:(A\ot B)\ot C\cong A\ot(B\ot C)$.
An $S$-morphism is a functor preserving the tensor product up to
coherent isomorphism. Finally we may consider the finitary 2-functor
$B:\Cat\to\Cat$ sending \C to $\C^4$. There are two $HB$-algebra
structures on an $S$-algebra \C, involving the derived operations
$((A\ot B)\ot C)\ot D$ and $A\ot(B\ot(C\ot D))$, and these induce
two monad maps $f',g':HB\to S$. The two isomorphisms $((A\ot B)\ot
C)\ot D\cong A\ot(B\ot(C\ot D))$ which can be built out of $\alpha$
induce two monad transformations $\phi,\psi:f'\to g'$, and we can
now form the {\em coequifier} of these 2-cells, namely the universal
monad map $q:S\to T$ for which $q\phi=q\psi$. A $T$-algebra is now
exactly what one might call a {\em semigroupoidal category}: a
category \C equipped with a tensor product $\ot$ which is
associative up to coherent isomorphism (coherent in the sense of the
Mac~Lane pentagon), but not necessarily having a unit. A
$T$-morphism is a strong semigroupoidal functor (one which preserves
the tensor product up to coherent isomorphism).

\pgph
Since \Mndf is monadic over \Endf via a finitary 2-monad $M$, we have
a lifted model structure on $\Mndf(=\malgs)$. As usual the cofibrant
objects are the flexible algebras, here called flexible monads.

We know that (i) free monads (on a finitary endo-2-functor) are
flexible, and that (ii) flexible colimits of flexible monads are
flexible. Since co-isoinserters and coequifiers are both flexible
colimits (a co-isoinserter can be constructed out of coinserters and
coequifiers) it follows that $T$ is a flexible monad. The key
feature of the presentation given above is that it used coinserters
and coequifiers but not such ``inflexible'' colimits as
coequalizers. As observed by Kelly, Power, and various
collaborators, a 2-monad is always flexible if it can be given by a
presentation which ``involves no equations between objects''. Thus
the 2-monad for monoidal categories is flexible, while that for
strict monoidal categories is not (it involves the equation
$A\ot(B\ot C)=(A\ot B)\ot C$.

\pgph
A monad morphism $k:S\to T$ induces a 2-functor $k^*_s:\talgs\to\salgs$
commuting with the forgetful 2-functors into \K; there is also an induced
2-functor $\talg\to\salg$, considered in Section~\ref{sect:ss} below. Explicitly,
$k^*_s$ sends a $T$-algebra $(A,a:TA\to A)$ to the $T$-algebra $(A,a')$,
where $a'$ is the composite of $a$ and $fA:SA\to TA$. Since
fibrations and trivial fibrations in the 2-categories of algebras are defined
as in \K, and $k^*_s$ commutes with the forgetful 2-functors, $k^*_s$ preserves
fibrations and trivial fibrations. It also has a left adjoint \cite{BKP},
and so is the right adjoint part of a Quillen adjunction. In Section~\ref{sect:ss}
we shall find conditions under which it is a Quillen equivalence.

\pgph
As well as the 2-category $\Mndf=\malgs$ of 2-monads and strict morphisms,
we can also consider the 2-category \malg of
finitary 2-monads on \K and {\em pseudomorphisms} of monads. Explicitly,
a pseudomorphism from $T$ to $S$ consists of a 2-natural transformation
$f:T\to S$ which preserves the unit and multiplication up to isomorphism;
these isomorphisms are required to satisfy coherence conditions formally
identical to those for a strong monoidal functor. One of the main
reason for considering pseudomorphisms is for dealing with
{\em pseudoalgebras}; this technique goes back to \cite{Kelly-LaxAlg}.
A pseudoalgebra for a 2-monad $T$ is an object
$A$ equipped with a morphism $a:TA\to A$ satisfying the usual algebra
axioms up to coherent isomorphism. This may be expressed by saying that
the 2-natural $\alpha:T\to\<A,A\>$ corresponding to $A$ is a
pseudomorphism of monads. Thus we have bijective correspondences
between pseudo $T$-algebra structures on $A$, pseudomorphisms
$T\to\<A,A\>$, strict morphisms $T'\to\<A,A\>$, and (strict)
$T'$-algebra structures on $A$, and in fact \ttalgs is isomorphic to
the 2-category \pstalgs of pseudo $T$-algebras and {\em strict} $T$-morphisms,
and similarly \ttalg isomorphic to the 2-category \pstalg of pseudo
$T$-algebras and pseudo $T$-morphisms.

\pgph\label{pgph:flexible-transport}
Notice that if $T$ is flexible, then by
Proposition~\ref{prop:flexibledomain} every pseudomorphism
$T\to\<A,A\>$ is isomorphic to a strict one. When this is translated
into a statement about algebras it states that for a flexible $T$,
every pseudoalgebra structure on an object $A$ is isomorphic in
\pstalg to a strict algebra structure on $A$ via pseudomorphism of
the form $(1_A,\phi)$. In particular, for a flexible monad, every
pseudoalgebra is isomorphic to a strict one.


\section{Structure and semantics}
\label{sect:ss}


\pgph We now turn from monads to their algebras. Whereas earlier we
considered model structures on \talgs, as a vehicle to understanding
the more important 2-category \talg, we now focus on \talg itself.
The passage from a 2-monad $T$ on \K to the 2-category \talg with
forgetful 2-functor $U:\talg\to\K$ is functorial. Given a (strict)
morphism $k:S\to T$ of 2-monads, the 2-functor $k^*_s:\talgs\to\salgs$
extends to a 2-functor $k^*:\talg\to\salg$, also commuting with the
forgetful 2-functors. In order to capture this situation, we consider the (enormous)
2-category \Twocat of (not necessarily small) 2-categories,
2-functors, and 2-natural transformations (ignoring the further
structure which makes it into a 3-category), and then the slice
2-category \slice, an object of which is a 2-category \LL equipped
with a 2-functor $U:\LL\to\K$, and a morphism of which is a
commutative triangle. If $M$ and $N$ are morphisms from $U:\LL\to\K$
to $U':\LL'\to\K$, a 2-cell from $M$ to $N$ is a 2-natural
transformation $\rho:M\to N$ whose composite with $U'$ is the
identity on $U$. Then there is a functor $\sem:\Mndf\op\to\slice$
which sends a 2-monad $T$ to $U:\talg\to\K$, and a morphism $j:S\to
T$ in \Mndf to $k^*:\talg\to\salg$.

\pgph Although size problems prevent there from being a model
structure on \slice, there are nonetheless obvious notions of
fibration and weak equivalence, which we now describe.

There is a Quillen model structure on the category
\twocat of small 2-categories and
2-functors, described in \cite{qm2cat}, for which the
weak equivalences are the biequivalences: these are the 2-functors
$F:\K\to\LL$ for which each functor $F:\K(A,B)\to\LL(FA,FB)$ is an
equivalence and furthermore for each $C\in\LL$ there is an $A\in\K$
and an equivalence $C\simeq FA$ in \LL. A 2-functor $F:\K\to\LL$ is
a fibration if each $F:\K(A,B)\to\LL(FA,FB)$ is a fibration in \Cat,
and moreover equivalences lift through $F$ in a sense made precise
in \cite{qmbicat}. (There is a mistake in the description of fibrations
and trivial cofibrations; this is corrected in \cite{qmbicat}, which
also provides a model structure on the category of bicategories and
strict homomorphisms, and shows that these two model categories are
Quillen equivalent.) The cofibrations are of course the maps with the left lifting
property with respect to the trivial fibrations; these are characterized
in \cite{qm2cat}.

Clearly the definitions of weak equivalence and fibration have nothing
to do with size, and one can easily extend them to give notions of
fibration and weak equivalence in \Twocat.

There is an evident functor $D:\slice\to\Twocat$ sending
$U:\LL\to\Cat$ to \LL, and we define a morphism $f$ in \slice to be a
weak equivalence or fibration if and only if $Df$ is one in \Twocat.

\pgph
Let $k:S\to T$ be a monad morphism. We consider the following induced
maps. First of all there is $k^*_s:\talgs\to\salgs$, which is a right
Quillen 2-functor. Then there is the 2-functor $k^*:\talg\to\salg$
which extends $k^*_s$. Finally there is the $\Ho\Cat$-functor
$\Ho(k^*_s):\Ho\talgs\to\Ho\salgs$. Since all objects in \talgs are
fibrant, $\Ho(k^*_s)$ is induced directly from $k^*_s$ without having
to use fibrant approximation. Thus $\Ho(k^*_s)$ is simply the underlying
$\Ho\Cat$-functor of $k^*:\talg\to\salg$.

\begin{proposition}
The following are equivalent:
\begin{enumerate}[(i)]
\item $k^*_s$ is a Quillen equivalence;
\item $\Ho(k^*_s)$ is an equivalence of $\Ho\Cat$-categories;
\item $k^*$ is a biequivalence.
\end{enumerate}
\end{proposition}

\proof
The equivalence of (i) and (ii) was proved in Section~\ref{pgph:Quillen-equivalence}.
For $T$-algebras $A$ and $B$, we have $k^*:\talg(A,B)\to\salg(k^*A,k^*B)$ an equivalence
if and only if $\Ho(k^*_s):\Ho\talgs(A,B)\to\Ho\salgs(k^*A,k^*B)$ is invertible,
while for a $T$-algebra $A$ and an $S$-algebra $C$, we have $k^* A\cong C$
in $\Ho\salgs$ if and only if $k^*A\simeq C$ in \salg. This gives the equivalence
between (ii) and (iii).
\endproof

\pgph Colimits in \Mndf of course become limits in $\Mndf\op$, and
cofibrations and weak equivalences in \Mndf become fibrations and
weak equivalences in $\Mndf\op$. Size issues notwithstanding, the
functor $\sem:\Mndf\op\to\slice$ sends limits to limits, fibrations
to fibrations, and trivial fibrations to trivial fibrations, as we
verify below, using the constructions $\<A,A\>$, $\{f,f\}$, and
$[\rho,\rho]$ of Section~\ref{pgph:<A,A>}.

\pgph To say that $\sem$ preserves limits is to say that it
sends colimits in \Mndf to limits in \slice. Consider the case of
coproducts. A product in \slice is just a fibre product in \twocat
(over \Cat). Suppose then that $(S_i)_{i\in I}$ is a small family of
finitary 2-monads on \Cat, with coproduct $S=\sum_i S_i$. The
product in \slice of the $\sem(S_i)$, is the 2-category in which an
object is a \K-object $A$ equipped with an $S_i$-algebra structure
$a_i:S_iA\to A$ for each $i$; a morphism between two such is a
\K-morphism $f:A\to B$ equipped with, for each $i$, an isomorphism
$\bar{f}_i:b_i.Sf\cong fa_i$ making $f$ into an $S_i$-morphism; and
a 2-cell between two such is a \K-transformation $f\to g$ compatible
with the $S_i$-morphism structure for each $i$. So to make $A$ into
an object of $\prod_i\sem(S_i)$ is to give a monad map
$\alpha_i:S_i\to\<A,A\>$ for each $i$; but this is precisely to give
a single monad map $\alpha:S\to\<A,A\>$; that is, an $S$-algebra
$a:SA\to A$ structure for $A$. The case of morphisms is treated
similarly. Let $(A,a)$ and $(B,b)$ be $S$-algebras, with notation
for the other associated maps as above. Then to make a \K-morphism
$f:A\to B$ into a morphism in $\prod_i\sem(S_i)$ is to give monad
maps $\phi_i:S_i\to\{f,f\}$ for each $i$, compatible with the
$\alpha_i:S_i\to\<A,A\>$ and $\beta_i:S_i\to\<B,B\>$. But by the
universal property of $S$, this amounts to giving a single monad map
$\phi:S\to\{f,f\}$ compatible with $\alpha$ and $\beta$; that is, to
a single $\bar{f}:b.Sf\cong fa$ making $f$ into an $S$-morphism.
Thus $\sem(S)$ and $\prod_i\sem(S_i)$ have the same objects and
morphisms; it remains to check the 2-cells, and this can be done
entirely analogously, using the construction $[\rho,\rho]$.

This proves that $\sem:\Mndf\op\to\slice$ preserves products; the
case of general limits is similar, and left to the reader.

\pgph Let $j:S\to T$ be a trivial cofibration in \Mndf, and so a
trivial fibration in $\Mndf\op$ from $T$ to $S$; we shall show that
$j^*:\talg\to\salg$ is a trivial fibration in \Twocat, and so that
$\sem(j)$ is a trivial fibration in \slice. Since $S$ (like every
other object of \Mndf) is fibrant, we know by
Proposition~\ref{prop:Cat-model} that
$\Mndf(j,S):\Mndf(T,S)\to\Mndf(S,S)$ is a surjective equivalence.
Thus there is a monad morphism $g:T\to S$ with $gj=1$, and a unique
isomorphism $\rho:jg\cong 1$ in \Mndf with $g\rho=\id$ and $\rho
j=\id$. By functoriality of \sem, we have $j^*g^*=1$ and
$g^*j^*\cong 1$, so $j^*$ is not just a trivial fibration, but in
fact a retract equivalence in \slice. This proves that
$\sem:\Mndf\op\to\slice$ sends trivial fibrations to trivial
fibrations.

\begin{remark}
In fact by the same argument even the trivial cofibrations for the
trivial model structure in \Mndf are sent to trivial fibrations in
\slice.
\end{remark}

\pgph The case of general weak equivalences is more delicate. A
morphism $f:S\to T$ in \Mndf is a weak equivalence if and only if
there exists a 2-natural $g:T\to S$ with $gf\cong 1$ and $fg\cong
1$. This $g$ will automatically be a monad pseudomorphism, but need
not in general be a monad morphism, thus it need not induce a 2-functor
$g^*:\salg\to\talg$. As a special case, consider the
weak equivalence $q:T'\to T$. Then $q^*$ is the inclusion
$\talg\to\pstalg$, which is a weak equivalence if and only if every
pseudo $T$-algebra is equivalent to a strict one: the ``general
coherence problem'' for $T$-algebras. This is still an open problem
in the current generality, but has been solved in various special
cases
--- see \cite{codescent} and the references therein.

If the monads $S$ and $T$ are flexible, however, then any weak equivalence
$f:S\to T$ does induce a biequivalence $f^*:\talg\to\salg$. This can be
proved using the observation of Section~\ref{pgph:flexible-transport}, or
it could also be deduced using Ken Brown's lemma \cite[1.1.2]{Hovey-book}.

On the other hand, the weak equivalences in \Mndf for the trivial
model structure are just the equivalences in \Mndf, and these {\em are}
mapped to weak equivalences in \slice.

\pgph The next thing to do is to check whether \sem preserves
fibrations; that is, whether $j^*:\talg\to\salg$ is a fibration in
\Twocat whenever $j:S\to T$ is a cofibration in \Mndf. First we
check that equivalences can be lifted through $j^*$. Suppose then
that $(A,a)$ is a $T$-algebra, that $(B,b)$ is an $S$-algebra, and
that $(f,\bar{f}):(B,b)\to j^*(A,a)$ is an equivalence of
$S$-algebras; the latter implies in particular that $(f,\bar{f})$ is
a morphism of $S$-algebras with $f:A\to B$ an equivalence in \K. We
must show that the $S$-algebra structure on $B$ can be extended to a
$T$-algebra structure in such a way that $f$ becomes a morphism of
$T$-algebras.

To do this, let $\beta:S\to\<B,B\>$ and $\alpha:T\to\<A,A\>$ be the
monad morphisms corresponding to $b:SB\to B$ and $a:TA\to A$. We
shall need, among other things, to extend $\beta$ along $j:S\to T$.
Let $\phi:S\to\{f,f\}$ be the monad morphism corresponding to
$(f,\bar{f})$. Then the diagram
$$\xymatrix{
S \ar[r]^{\phi} \ar[d]_{m} & \{f,f\} \ar[d]^{d} \\
T \ar[r]_{\alpha} & \<A,A\> }$$
of monad morphisms commutes, and the equivalence lifting property for
$j^*:\talg\to\salg$ now amounts to the existence of a fill-in. Since
$j$ was assumed to be a cofibration, this fill-in will exist provided
that $f:\{f,f\}\to\<A,A\>$ is a trivial fibration in \Mndf, or equivalently
in \Endf. But $f:A\to B$ was assumed an equivalence, thus $\<f,B\>$ is
an equivalence, and the result now follows by the general fact that if
$$\xymatrix @R1pc {
W \ar[r] \ar[dd]_{p} & X \ar[dd]^{w} \\
{}\rtwocell<\omit>{} & {} \\
Y \ar[r] & Z }$$
is an iso-comma object (in any 2-category) and $w$ an equivalence then
$p$ is a retract equivalence.

This proves the equivalence lifting property. We now turn to the 2-dimensional
aspect. This asserts that if $(f,\bar{f}):(A,a)\to(B,b)$ is a $T$-morphism,
and $\rho:(g,\bar{g})\cong j^*(f,\bar{f})$ is an invertible $S$-transformation,
then we can lift this to an invertible $T$-transformation. This is
entirely straightforward and is true for {\em any} monad morphism $j:S\to T$.

This completes the verification that $\sem:\Mndf\op\to\slice$
preserves fibrations and trivial fibrations.

\pgph
The statement about fibrations is that if $j:S\to T$ is a cofibration
in \Mndf, then $j^*:\talg\to\salg$ is a fibration in \Twocat. As a special
case, if $T$ is flexible (cofibrant), then the forgetful 2-functor
$U:\talg\to\K$ is a fibration in \Twocat. This amounts to the facts that (i)
if $(A,a)$ is a $T$-algebra, and $f:B\to A$ an equivalence in \K, then there
is a $T$-algebra structure $(B,b)$ on $B$, for which $f$ can be made into
a $T$-morphism, and (ii) if $(f,\bar{f}):(A,a)\to(B,b)$ is a $T$-morphism,
and $\phi:g\cong f$, then $g$ can be made into a $T$-morphism $(g,\bar{g})$
in such a way that $\phi$ is a $T$-transformation $(g,\bar{g})\to(f,\bar{f})$.
As before, (ii) is true for any 2-monad $T$, while (i) asserts that $T$-algebra
structure can be transported along equivalences. This is true for flexible
monads, but not in general. It is false, for example, in the case of the 2-monad
$T$ for strict monoidal categories.
Consider the category $C$ of countable sets, viewed as a monoidal category
under the cartesian product. This can be replaced by an equivalent strict
monoidal category $A$. If $B$ be a skeleton of $C$ (choose one set of each
countable cardinality), then there exists an equivalence of categories $f:B\to A$.
But by an argument due to Isbell (see\cite[VII.1]{CWM}), there is no way to transport
the strict monoidal structure on $A$ across the equivalence $f$ to obtain a
strict monoidal structure on $B$.

\pgph
It is only the hugeness of \slice that prevents it from having a model structure,
and \sem from having a left adjoint, and so becoming a right Quillen 2-functor.
It would be interesting to find a full sub-2-category of \slice containing the
image of \sem, admitting a model structure with the fibrations and weak equivalences
defined as in \slice, and on which a left adjoint to \sem can be defined.

\bibliographystyle{plain}

\end{document}